\documentclass{amsart}

\usepackage{amsmath, amssymb, enumerate}
\input xy
\xyoption{all}
\numberwithin{equation}{section}

\begin{document}

\newtheorem{theorem}{Theorem}[section]
\newtheorem{lemma}[theorem]{Lemma}
\newtheorem{corollary}[theorem]{Corollary}
\newtheorem{fact}[theorem]{Fact}
\newtheorem{proposition}[theorem]{Proposition}

\theoremstyle{definition}
\newtheorem{example}[theorem]{Example}
\newtheorem{definition}[theorem]{Definition}
\newtheorem{question}[theorem]{Question}

\theoremstyle{remark}
\newtheorem{claim}[theorem]{Claim}
\newtheorem{remark}[theorem]{Remark}

\def\endo{\operatorname{End}}
\def\fix{\operatorname{Fix}}
\def\ccma{\operatorname{CCMA}}
\def\dcf{\operatorname{DCF}}
\def\acf{\operatorname{ACF}}
\def\dcfa{\operatorname{DCFA}}
\def\acfa{\operatorname{ACFA}}
\def\ccm{\operatorname{CCM}}
\def\dccm{\operatorname{DCCM}}
\def\jet{\operatorname{Jet}}
\def\hom{\operatorname{Hom}}
\def\Th{\operatorname{Th}}
\def\cb{\operatorname{Cb}}
\def\cbsigma{\operatorname{Cb}_{\sigma}}
\def\qfcb{\operatorname{qfCb}}
\def\tp{\operatorname{tp}}
\def\dimnabla{\operatorname{dim}_{\nabla}}
\def\qftp{\operatorname{qftp}}
\def\Stab{\operatorname{Stab}}
\def\stp{\operatorname{stp}}
\def\acl{\operatorname{acl}}
\def\aclnabla{\operatorname{acl}_{\nabla}}
\def\dcl{\operatorname{dcl}}
\def\dclnabla{\operatorname{dcl}_{\nabla}}
\def\eq{\operatorname{eq}}
\def\loc{\operatorname{CCM-loc}}
\def\jet{\operatorname{Jet}}
\def\alg{\operatorname{alg}}
\def\rank{\operatorname{rank}}
\def\SU{\operatorname{SU}}
\def\aut{\operatorname{Aut}}
\def\gal{\operatorname{Gal}}
\def\loc{\operatorname{loc}}
\def\pr{\operatorname{pr}}
\def\ga{\operatorname{GA}}

\def\ccmforallnabla{\ccm_{\forall,\nabla}}

\def\bX{\overline X}

\def\C{\mathcal C}
\def\A{\mathcal A}
\def\M{\mathcal M}
\def\N{\mathcal N}
\def\F{\mathcal F}
\def\O{\mathcal O}
\def\P{\mathcal P}
\def\U{\mathcal U}
\def\K{\mathcal K}

\renewcommand{\L}{{\mathcal L}}

\def\PP{{\mathbb P}}
\def\UU{\mathbb{U}}
\def\QQ{\mathbb{Q}}
\def\ZZ{\mathbb{Z}}
\def\NN{\mathbb{N}}
\def\CC{\mathbb{C}}

\newcommand{\spanof}[1]{\left< #1\right>}

\newcommand{\slantfrac}[2]{\hbox{$\,^{#1}\!/_{#2}$}}  
\newcommand{\quot}[2]{\slantfrac{#1}{#2}}

\newcommand{\isom}{\cong}
\newcommand{\elemequiv}{\equiv}
\newcommand{\elemeq}{\elemequiv}

\newcommand{\maps}{\rightarrow} 

\newcommand{\imp}{\rightarrow}
\providecommand{\s}{\sigma}
\providecommand{\id}{\operatorname{id}}

\renewcommand{\O}{\mathcal{O}}

\newcommand{\elex}{\preccurlyeq}
\newcommand{\elres}{\succcurlyeq}


\def\Ind#1#2{#1\setbox0=\hbox{$#1x$}\kern\wd0\hbox to 0pt{\hss$#1\mid$\hss}
\lower.9\ht0\hbox to 0pt{\hss$#1\smile$\hss}\kern\wd0}
\def\ind{\mathop{\mathpalette\Ind{}}}
\def\notind#1#2{#1\setbox0=\hbox{$#1x$}\kern\wd0\hbox to 0pt{\mathchardef
\nn=12854\hss$#1\nn$\kern1.4\wd0\hss}\hbox to
0pt{\hss$#1\mid$\hss}\lower.9\ht0 \hbox to
0pt{\hss$#1\smile$\hss}\kern\wd0}
\def\nind{\mathop{\mathpalette\notind{}}}

\date{\today}

\title[Model theory of meromorphic vector fields]{A model theory for meromorphic vector fields}
\author{Rahim Moosa}

\address{Rahim Moosa\\
University of Waterloo\\
Department of Pure Mathematics\\
200 University Avenue West\\
Waterloo, Ontario \  N2L 3G1\\
Canada}
\email{rmoosa@uwaterloo.ca}

\thanks{R. Moosa was partially supported by an NSERC Discovery Grant.}



\begin{abstract}
Motivated by the study of meromorphic vector fields, a model theory of ``compact complex manifolds equipped with a generic derivation" is here proposed.
This is made precise by the notion of a  {\em differential $\ccm$-structure}.
A first-order axiomatisation of existentially closed differential $\ccm$-structures is given.
The resulting theory, $\dccm$, is a common expansion of the theories of differentially closed fields and compact complex manifolds.
A study of the basic model theory of $\dccm$ is initiated, including proofs of completeness, quantifier elimination, elimination of imaginaries, and total transcendentality.
The finite-dimensional types in $\dccm$ are shown to be precisely the generic types of meromorphic vector fields.
\end{abstract}

\maketitle

\section{Introduction}

\noindent
The model-theoretic approach to systems of (ordinary) algebraic differential equations is via the first-order theory of differentially closed fields in characteristic zero ($\dcf_0$).
Such systems of equations, at least in the autonomous case when the equations have constant parameters, can be presented geometrically as algebraic vector fields; namely a projective algebraic variety $X$ equipped with a rational section $v:X\to TX$ to the tangent space.
In fact, the finite-dimensional fragment of $\dcf_0$ essentially coincides with the birational geometry of algebraic vector fields.
(See, for example, \cite{6lectures} for an exposition of $\dcf_0$ from this point of view.)
Here, I am interested in generalising this model-theoretic framework to {\em meromorphic} vector fields; namely when $X$ is a compact complex manifold that is not necessarily algebraic and $v$ is a meromorphic section to the holomorphic tangent bundle.
While $\dcf_0$ is built on the theory of algebraically closed fields ($\acf_0$), the new theory I am seeking should be built on a first-order theory of compact complex manifolds.

About thirty years ago, in~\cite{easter}, as part of the development of the notion of ``Zariski-type structure", Zilber proposed a model theory for compact complex manifolds.
Unlike $\acf_0$ and $\dcf_0$, the first-order theory proposed by Zilber for compact complex manifolds was not given by an explicit axiomatisation, nor as the model companion of a natural class of algebraic structures, but rather as theories of particular structures:
a compact complex manifold~$M$ is viewed as a first-order structure in the language where there is a predicate for each closed complex-analytic subset of each finite cartesian power of $M$.
Zilber showed that the theory of any such structure shares many properties with its algebraic predecessors: in particular, they admit quantifier elimination and are of finite Morley rank (bounded by the dimension of $M$).
Later, in the work of Hrushovski~\cite{udi-icm} and Pillay~\cite{pillay00}, for example, it became common to consider all compact complex manifolds, indeed all (reduced) compact complex-analytic spaces, at once, in a multisorted structure whose theory now goes by the name $\ccm$.
Like differentially closed fields, $\ccm$ is a proper expansion of $\acf_0$.
Also like $\dcf_0$, much of the richness of geometric stability theory absent in $\acf_0$ is present in $\ccm$.
For example, all cases of the Zilber trichotomy appear.

In this paper, I will present a common expansion of $\ccm$ and $\dcf_0$, which I will call $\dccm$.
It will turn out (in Section~\ref{fd}, below) that the finite-dimensional fragment of $\dccm$ will capture, precisely, the bimeromorphic geometry of meromorphic vector fields.
As such, it achieves the goal set out in this introduction.

The theory $\dccm$ arises by considering {\em differential $\ccm$-structures}, essentially by adding a ``derivation" to the definable closure of a generic point of a sort, say $X$, in $\ccm$.
This can be made sense of because the elements of the definable closure of a generic point of $X$ can be viewed as meromorphic maps from $X$ to other sorts, and hence can be differentiated.
See Sections~\ref{mvts} and ~\ref{ds} for a detailed explanation.

The specific goals of this paper are:
\begin{enumerate}
\item
to show that the (universal) theory of differential $\ccm$-structures admits a model companion, which is $\dccm$, by giving a geometric first-order axiomatisation of the existentially closed models (Theorem~\ref{dccm});
\item
to show that $\dccm$ is complete, admits quantifier elimination (Proposition~\ref{apqe}) and elimination of imagainaries (Theorem~\ref{ei}), and to give a geometric characterisation of definable and algebraic closure (Proposition~\ref{dclacl});
\item
to show that $\dccm$ is totally transcendental (Theorem~\ref{lambdastable}), and to give geometric characterisations of nonforking independence (Corollary~\ref{stable}); and,
\item
to establish the correspondence between finite-dimensional types (over the empty set) in $\dccm$ and meromorphic vector fields (Theorem~\ref{fdmervf}).
\end{enumerate}
The proofs proceed largely by finding geometric analogues for the algebraic arguments already familiar from $\dcf_0$.

The next step in the study of $\dccm$, not attempted here, would be to establish the canonical base property for finite dimensional types following the strategy of~\cite{pillayziegler03} in the case of $\dcf_0$.
This will involve developing a theory of jet spaces in $\dccm$, see for example~\cite{ccma} where this was done for compact complex manifolds with a generic automorphism ($\ccma$).
In any case, once the canonical base property is established, a concrete manifestation of the Zilber dichotomy for finite dimensional minimal types in $\dccm$ will follow.
It would then be reasonable to expect that many of the recent applications of model theory to algebraic vector fields, as carried out in~\cite{C3C2} and~\cite{abred} for example, would extend to meromorphic vector fields.

The process of adding an automorphism to any given first-order theory of interest, and then seeking a model companion, is well-studied (see~\cite{chatzidakis-pillay}).
Here we have ``added a derivation" instead.
Clearly, this does not make sense for an arbitrary theory.
But following the ideas presented here, it may be worth investigating a robust general setting where adding a derivation does make sense.
A likely candidate might be that of Zariski-type structures in Zilber's sense; one that expands $\acf_0$ and admits a functor that extends to all sorts the tangent space construction on algebraic varieties.
\vfill
\pagebreak

\section{Meromorphic varieties and their tangent spaces}
\label{mvts}

\noindent
In this section I want to slightly loosen the usual formalism for doing the model theory of compact complex-analytic spaces.
The idea is to allow us to work directly in the ``compactifiable" rather than compact setting.

For the fundamental notions from complex-analytic geometry we suggest~\cite{fischer76}.
Given a reduced compact complex-anlaytic space $X$, by the {\em Zariski topology} on~$X$ we will mean the (noetherian) topology of closed complex-analytic subsets  of~$X$.
This will not conflict with the usual meaning of the Zariski topology in the case that $X$ is a projective complex-algebraic variety, because in that case the complex-analytic and complex-algebraic sets agree (Chow's Theorem).

\begin{definition}
By a {\em meromorphic variety}, $X$,  we will mean a Zariski dense and open subset of a reduced compact complex analytic space which we will tend to denote by $\overline X$.
Of course, $X$ inherits from $\bX$ the structure of a reduced complex-analytic space in its own right, and may admit other (bimeromorphically equivalent)  compactifications.
But we view $X$ as embedded in the given compactification $\overline X$.

Cartesian products of meromorphic varieties, $X\times Y$, are viewed as meromorphic varieties with the compactification $\overline{X\times Y}=\bX\times\overline Y$.

By the {\em Zariski topology} on $X$ we mean the topology induced by the Zariski topology on $\bX$.
Note that this is a coarser topology than that of the closed complex-analytic subsets of $X$; such a set is Zariski closed in $X$ if and only if its (euclidean) closure in $\bX$ is Zariski closed.

By a {\em definable holomorphic map} $f:X\to Y$ of meromorphic varieties we mean a holomorphic map that extends to a meromorphic map $\overline f:\bX\to\overline Y$.
Equivalently, the graph of $f$ is Zariski closed in $X\times Y$.
More generally, a {\em definable meromorphic map} $f:X\to Y$ is a meromorphic map that extends to a meromorphic map from $\bX$ to $\overline Y$.
Such a map is {\em dominant} if its image is Zariski dense in $Y$.
\end{definition}

The usual model-theoretic set-up is to consider the first-order theory of the multi-sorted structure $\A$ where there is a sort for each reduced and irreducible compact complex-analytic space, and a predicate for each Zariski closed subset of each finite cartesian product of sorts.
See, for example, the surveys \cite{survey} and~\cite{kaehler-survey}.
Every meromorphic variety, in the above sense, is $0$-definable in $\A$, as is every Zariski closed subset of every finite cartesian product of meromorphic varieties.
It therefore does no harm to work instead with the expansion of $\A$ to the multi-sorted structure $\M$ where there is a sort for each irreducible meromorphic variety, and a predicate for each Zariski closed subset of each finite cartesian product of sorts.
So we have added some sorts and some predicates, but they were all already $0$-definable in the original structure.
I will denote by $L$ the language of $\M$, and by $\ccm$ the first-order $L$-theory of $\M$.
It admits quantifier elimination and elimination of imaginaries, and, sort by sort, is of finite Morley rank.

Every quasi-projective complex-algebraic variety $V$, given with an embedding in a projective compactification $\overline V$, is a meromorphic variety, and the algebraic and analytic Zariski topologies on $V$ agree.
In particular, definable holomorphic maps in this case are just regular morphisms, and definable meromorphic maps are rational.
In this way, algebraic geometry lives as a pure reduct of $\ccm$.

Our main use of the flexibility that $\M$ affords is that the collection of sorts is closed under taking tangent spaces.
Recall that the {\em tangent space} of a complex-analytic space $X$ is the linear fibre space $\pi:TX\to X$ associated to the sheaf of differentials, $\underline\Omega_X$, on $X$.
So $TX$ is a complex-analytic space and $\pi:TX\to X$ is a surjective holomorphic map whose fibres are uniformly equipped with the structure of a complex vector space, in the sense that there are holomorphic maps for addition, $+:TX\times_XTX\to TX$, scalar multiplication $\lambda:\mathbb C\times TX\to TX$, and zero section $z:X\to TX$, all over $X$, satisfying the vector space axioms.
For any point $p\in X$, the {\em tangent space to $X$ at $p$} is the fibre of $\pi:TX\to X$ above $p$, denoted by $T_pX$, and it is canonically isomorphic as a complex vector space to $\operatorname{Hom}_{\mathbb C}(\mathfrak m_{X,p}/\mathfrak m_{X,p}^2,\mathbb C)$, where $\mathfrak m_{X,p}$ is the maximal ideal of the local ring of $X$ at $p$.

We claim that when $X$ is a meromorphic variety so is $TX$, and that  $\pi, +,\lambda, z$ are all definable holomorphic maps.
Let us first consider the case when $X=\bX$ is already compact.
We are looking for a natural compactification of $TX$.
In fact, there is a canonical way to do this for any linear fibre space $\L(\F)\to X$ associated to a coherent analytic sheaf $\F$ on $X$; it is just the relativisation of the usual embedding of $\CC^n$ in the projectivisation of $\mathbb C^{n+1}$.
One considers the coherent analytic sheaf $\F\times\O_X$ of rank one greater than $\F$, and then the associated projective linear space $\PP(\F\times\O_X)\to X$.
See~\cite[Section~1.9]{fischer76} for details.
Then $\L(\F)$ embeds in $\PP(\F\times\O_X)$ over $X$ as a Zariski open set in such a way that the linear structure (namely, $\pi, +,\lambda, z$) extends meromorphically to the projective linear space.
Applying this to $\F=\underline\Omega_X$ gives $TX\to X$ the meromorphic structure we are looking for, namely $\overline{TX}:=\PP(\underline\Omega_X\times\O_X)$.
Now, if we consider a general meromorphic variety $X$ embedded in $\bX$, then the linear space $TX\to X$ is just the restriction to $X$ of $T\bX\to \bX$, hence $\overline{T\bX}$ will serve as a compactification for $TX$, to which the linear structure extends meromorphically.

Recall that the tangent space construction is functorial:
For each meromorphic (respectively, holomorphic) map $g:X\to Y$ between complex-analytic spaces there is a meromorphic (respectively, holomorphic) map $dg:TX\to TY$, such that
$$\xymatrix{
TX\ar[r]^{dg}\ar[d]_{\pi_X}&TY\ar[d]^{\pi_Y}\\
X\ar[r]^g&Y
}$$
commutes, and we have the functoriality property $d(g\circ h)=(dg)\circ(dh)$.
If $X$ and $Y$ are meromorphic varieties, and $g:X\to Y$ is definable meromorphic (respectively, holomorphic), then so is $dg:TX\to TY$.
That is, if $g$ extends to a meromorphic map $\bX\to\overline Y$ the $dg$ extends to a meormorphic map $\overline{TX}\to\overline{TY}$.

\bigskip
\section{The differential structure}
\label{ds}
\noindent
By a {\em $\ccm$-structure} I will mean a definably closed subset of a model of $\ccm$.
In other words, a model of $\ccm_{\forall}$.
The goal of this section is to describe what we might consider a ``derivation" on a $\ccm$-structure.
But first, let us recall what $\ccm$-structures themselves look like.

Since we are in a relational language in which all elements of $\M$ are named, a  model of $\ccm_\forall$ is simply a nonempty subset $A$ of an elementary extension $\N$ of~$\M$.
As we are in a multi-sorted setting, nonempty is meant relative to every sort: so $S(A)\neq\emptyset$ for all sorts $S$ of $L$.
But we will mostly be interested in finitely generated definably closed substructures, so where $A=\dcl(a)$ for some $a\in X(\N)$ and some irreducible meromorphic variety~$X$.
Replacing $X$ by the locus of~$a$, we may assume that $a$ is a generic point of~$X$ in the sense that it is not contained in $Y(\N)$ for any proper Zariski closed subset $Y\subsetneq X$.
In that case we can identify $A$ with the set of all definable meromorphic maps $g:X\to S$ as $S$ ranges over all other sorts.
Indeed the identification is given by $g\mapsto g(a)\in S(A)$, noting that every point of $S(A)$ arises this way as $A=\dcl(a)$, and that if two definable meromorphic maps agree on $a$ then they agree on $X$ by genericity.

It is worth comparing to the algebraic case, so when $X$ happens to be a quasi-projective complex-algebraic variety. In that case one only needs to consider the single target sort $S=\PP$, the projective line.
Indeed, in that case, $\dcl(a)=\CC(X)$ is just the field of rational functions.
For nonalgebraic meromorphic varieties, if we only considered $S=\PP$ we would obtain the {\em meromorphic function field} of $X$, and not necesarilly the full definable closure of a generic point.
Indeed, on some compact complex-analytic spaces, namely those of {\em algebraic dimension} $0$, there are no nonconstant meromorphic functions, but many nonconstant meromorphic maps to other sorts.

The differential structure I want to consider is motivated by the study of the following natural objects in bimeromorphic geometry:

\begin{definition}
\label{mervf}
By a {\em meromorphic vector field} we will mean an irreducible meromorphic variety $X$ equipped with a definable meromorphic section $v:X\to TX$ to the tangent space of $X$.
\end{definition}

For example, every meromorphic variety equipped with its zero section is a meromorphic vector field, which we call the {\em trivial vector field}.
On the other hand every rational vector field on an irreducible quasi-projective complex-algebraic variety is a meromorphic vector field.
For nontrivial, nonalgebraic examples, one can start with any compact complex manifold $X$ and consider the connected component of its automorphism group $G=\aut_0(X)$, which is itself a complex Lie group.
The Lie algebra of $G$ consists of holomorphic vector fields on $X$, see~\cite[III.1]{kobayashi}.
These give rise to many examples of nonalgebraic meromorphic vector fields.

Suppose $(X,v)$ is a meromorphic vector field, $\N\succeq\M$ is an elementary extension, and $a\in X(\N)$ is a generic point of $X$.
What structure does $v$ induce on $A:=\dcl(a)$?
Well, for any definable meromorphic $g:X\to S$, we have the definable meromorphic map $\nabla_v(g):=dg\circ v:X\to TS$.
Viewing $g\in S(A)$ we have defined a function $\nabla_v:S(A)\to TS(A)$, for all sorts $S$.
Here are two salient properties of this function that are easily verified using the functoriality of the tangent space construction:
\begin{itemize}
\item
$\pi\circ\nabla_v(g)=g$ where $\pi:TS(A)\to S(A)$ is the projection.
\item
$df\circ\nabla_v(g)=\nabla_v(f\circ g)$ for any definable meromorphic $f:S\to T$.\end{itemize}
We are thus lead to consider the following notion:

\begin{definition}
\label{ccmforallnabla}
Let $L_{\nabla}=L\cup\{\nabla\}$ where $\nabla=(\nabla_S:S\text{ sort of }L)$ and $\nabla_S$ is a function symbol from the sort $S$ to the sort $TS$.
Let $\ccm_{\forall,\nabla}$ denote the universal $L_\nabla$-theory which is obtained by adding to $\ccm_\forall$ the following axioms:
\begin{enumerate}
\item
For each sort $S$, $\nabla_S:S\to TS$ is a section to $\pi:TS\to S$.
\item
For each definable meromorphic map $f:S_1\to S_2$ between sorts, the following diagram commutes:
$$\xymatrix{
TS_1\ar[r]^{df}&TS_2\\
S_1\ar[u]^{\nabla_{S_1}}\ar[r]^f&S_2\ar[u]_{\nabla_{S_2}}
}$$
Remembering that $f$ and $df$ are not function symbols in the language but rather their graphs are predicates, what we mean by this is the axiom
$$\forall xy\big ((x,y)\in \Gamma(f)\implies(\nabla_{S_1} x,\nabla_{S_2} y)\in \Gamma(df)\big).$$
\end{enumerate}
We will usually drop the subscript and write $\nabla$ for $\nabla_S$ whenever it is clear from context which sort we are working in.
\end{definition}

One consequence of Axiom~(2) that gets used often without mention is that $\nabla(a_1,a_2)=(\nabla a_1,\nabla a_2)$ under the identification $T(S_1\times S_2)=TS_1\times TS_2$.

We can always extend uniquely to the definable closure:

\begin{proposition}
\label{extdcl}
Suppose $A\subseteq\N\models\ccm$ and $(A,\nabla)\models\ccmforallnabla$.
Then there is a unique extension of $\nabla$ to $\dcl(A)$ making it a model of $\ccmforallnabla$.
\end{proposition}

\begin{proof}
Let $B:=\dcl(A)$.
Given a sort $S$ we need to define $\nabla$ on $S(B)$.
Fix $b\in S(B)$ and let $X:=\loc(b)\subseteq S$ so that $b\in X(B)$ is generic.
Since $b\in\dcl(A)$, there exists some other irreducible meromorphic variety $Y$ admitting a dominant definable meromorphic map $f:Y\to X$, and a generic point $a\in Y(A)$, such that $b=f(a)$.
Now, $df:TY\to TX$ and $\nabla(a)\in TY(A)$.
Define $\nabla(b):=df(\nabla a)$.
Indeed, this is forced upon us by Axiom~(2) of Definition~\ref{ccmforallnabla}, and hence takes care of the uniqueness part of the statement.

We have to check that it is well-defined.
Suppose  we have another $f':Y'\to X$ and $a'\in Y'(A)$ generic such that $b=f'(a')$ as well.
Let $Z=\loc(a,a')\subseteq Y\times Y'$ and consider $\bar f:=(f,f'):Z\to X^2$.
Since $\bar f$ takjes a generic point of $Z$ to the diagonal $D\subseteq X^2$ we have that $\bar f(Z)\subseteq D$.
Hence $d\bar f:TZ\to T(X^2)$ lands in $TD$ which is the diagonal in $T(X^2)=(TX)^2$.
Since
$d\bar f(\nabla(a,a'))=(df(\nabla a),df'(\nabla a'))$ by functoriality, this means that $df(\nabla a)=df'(\nabla a')$, as desired.

Next, observe that $\nabla$ so defined is a function from $S(B)$ to $TS(B)$, and is a section to $\pi:TS\to S$.
That is, $(B,\nabla)$ does satisfy Axiom~(1) of Definition~\ref{ccmforallnabla}.
Taking $f=\id$ in the above construction, we see also that $(A,\nabla)\subseteq(B,\nabla)$.

It remains to verify Axiom~(2).
That is, given $g:S_1\to S_2$ a definable meromorphic map between sorts, and $b_i\in S_i(B)$ with $g(b_1)=b_2$, we need to show that $dg(\nabla b_1)=\nabla b_2$.
Note that by concatenating  -- namely working in cartesian products -- we can arrange things so that $b_1$ and $b_2$ are defined over the same tuple from~$A$.
That is, there is a sort $S$ with $a\in S(A)$ such that $b_1=f_1(a)$ and $b_2=f_2(a)$ where $f_i:S\to S_i$ are definable meromorphic maps.
Taking Zariski loci we may assume that $a$ is generic in $S$ and that each $b_i$ is generic in $S_i$.
Hence
\begin{eqnarray*}
dg(\nabla b_1)
&=& dg(df_1(\nabla a))\ \text{ by how $\nabla$ is defined on $B$}\\
&=& d(gf_1)(\nabla a))\ \text{ by functoriality}\\
&=&df_2(\nabla a))\ \ \ \text{ as $gf_1=f_2$, as that is the case on the generic $a$}\\
&=& \nabla b_2\ \ \ \ \text{ by how $\nabla$ is defined on $B$,}
\end{eqnarray*}
as desired.
\end{proof}

\begin{definition}
A {\em differential $\ccm$-structure} is a model $(A,\nabla)\models\ccmforallnabla$ such that $A=\dcl(A)$.
\end{definition}

As a consequence of Proposition~\ref{extdcl}, when working with models of $\ccmforallnabla$ there is little loss of generality in assuming that we have a differential $\ccm$-structure, namely that the underlying set is definably closed in $\ccm$.

It is worth observing that standard points are always constant:

\begin{lemma}
\label{zerom}
Suppose $(A,\nabla)$ is a differential $\ccm$-structure and $S$ is a sort.
If $p\in S(\M)$ then $\nabla(p)=0\in T_pX$.
\end{lemma}

\begin{proof}
Note that $X:=\{p\}$ is itself an irreducible meromorphic variety, and we can consider the containment as a definable holomorphic map $f:X\to S$.
Now $TX=\{(p,0)\}$, and hence $\nabla_X=0$.
But, by Axiom~(2) of Definition~\ref{ccmforallnabla}, this forces $\nabla_S(p)=df(\nabla_X(p))=0$ as $df_p:T_pX\to T_{f(p)}S$ is a linear map.
\end{proof}

In the finitely $\dcl$-generated case we recover precisely the meromorphic vector fields that motivated Definition~\ref{ccmforallnabla}:

\begin{proposition}
Suppose $X$ is an irreducible meromorphic variety, $a$ is a generic point of~$X$ in some elementary extension, and $A=\dcl(a)$.
Then the differential $\ccm$-structures on $A$ are precisely the $\nabla_v$ induced by meromorphic vector fields $v:X\to TX$.
\end{proposition}

\begin{proof}
We have already seen that $(A,\nabla_v)\models\ccm_{\forall,\nabla}$ if $(X,v)$ is a meromorphic vector field.
For the converse, suppose $(A,\nabla)\models\ccm_{\forall,\nabla}$.
Note that $a\in X(A)$ and $\nabla(a)\in TX(A)$.
As definable meromorphic maps, $a\in X(A)$ is the identitiy map on $X$ and $\nabla(a)\in TX(A)$ is some $v:X\to TX$.
Axiom~(1) ensures that $v$ is a section to $\pi:TX\to X$, and hence a meromorphic vector field on $X$.
It remains to verify that $\nabla=\nabla_v$.
Let $g(a)\in S(A)$ where $g:X\to S$ is a definable meromorphic map and $S$ is a sort.
Then
$$\nabla_v(g(a))=dg\circ v(a)=dg\circ\nabla(a)=\nabla(g(a))$$
where the final equality is by Axiom~(2).
\end{proof}

So the study of meromorphic vector fields amounts to the study of (finitely generated) differential $\ccm$-structures.
In the usual model-theoretic way, we will eventually look for a model companion, a theory that axiomatises the existentially closed differential $\ccm$-structures.

We conclude this section by extending the notion of differential $\ccm$-structure to a setting where $\nabla$ is allowed to take values in an extension.
This will be useful in what follows.

\begin{definition}
\label{ccmforallnablaext}
Suppose $\N\models\ccm$ and $A\subseteq\N$ is a definably closed set.
By an {\em $\N$-valued differential $\ccm$-structure} on $A$ we mean a map $\nabla:S(A)\to TS(\N)$, for every sort $S$, such that $\nabla$ is a section to $\pi:TS\to S$, and such that $df(\nabla a)=\nabla(f(a))$ for all $a\in S(A)$ and all definable meromorphic maps~$f$.
\end{definition}

\bigskip
\section{Prolongations}
\label{prolongations}

\noindent
In this section we construct a version of the tangent space that is twisted by a differential structure.
Since differential structure only has content in proper elementary extensions of $\M$, this will necessarily be about ``meromorphic varieties over parameters" in arbitrary models of $\ccm$, which we begin by reviewing.

Fix a model $\N\models\ccm$.
Given an irreducible meromorphic variety $X$, we view it as a sort of $L$ consider its $\N$-points $X(\N)$.
Let us recall the Zariski topology on $X(\N)$ with parameters from $\N$, sometimes referred to as the {\em nonstandard} Zariski topology to emphasise that we are not necessarily in the prime model $\M$.
See~\cite[Section~2]{ret} for a more detailed discussion.
Every Zariski closed subset $Y\subseteq X(\M)$ is named as a predicate in $L$ and so we can consider $Y(\N)$.
These are the $0$-definable Zariski closed subsets of $X(\N)$.
More generally, given a set of parameter $A\subseteq\N$,
a Zariski closed subset $Y\subseteq X(\N)$ over $A$ is a subset of the form $Y=Z_a$ where $a\in S(A)$ is a generic point of another sort $S$ and $Z\subseteq S\times X$ is a ($0$-definable) Zariski closed subset that projects dominantly on $S$.
In diagrams:
$$\xymatrix{
Z\ar@{^{(}->}[r]\ar[d]_{\rho} &S\times X\\
S}$$
That is, $Y=Z_a$ arises as the generic member of a $0$-definable family of Zariski closed subsets of $X$.
This forms a noetherian topology on $X(\N)$.
If $Y$ is $A$-irreducible then we can take $Z$ to be irreducible, and if $Y$ is absolutely irreducible the we can take $Z$ so that $\rho:Z\to S$ is a fibre space, meaning its general fibres in the standard model are irreducible.

The general standard fibres of $\rho: Z\to S$ will be of constant dimension when~$Z$ is irreducible, giving rise to a notion of dimension for irreducible Zariski closed subsets of $X(\N)$, which we denote by $\dim Y$.

The tangent space construction extends to nonstandard Zariski closed sets.
Fix $Y=Z_a$ as above.
Then the tangent spaces of the fibres of $\rho$ in the standard model vary uniformly: Consider the following diagram:
$$\xymatrix{
S\times X &Z\ar@{_{(}->}[l]\ar[d]_{\rho}& TZ\ar[l]\ar[d]^{d\rho}\ar@{^{(}->}[r] &TS\times TX\\
& S & TS\ar[l]
}$$
and let $z:S\to TS$ be the zero section.
For any $p\in S(\M)$, the fibre $(TZ)_{z(p)}\subseteq TX$ of $d\rho$ above $z(p)$ is nothing other than $T(Z_p)$, the tangent space of $Z_p\subseteq X$.
Hence we define the {\em tangent space of $Y=Z_a$} in $\N$, denoted $TY$, to be $(TZ)_{z(a)}$.

Suppose, now, that $A=\dcl(A)$ and we have a $\N$-valued differential $\ccm$-structure $\nabla$ on $A$.
Then, instead of considering the zero section, we can consider the differential section $\nabla$.
That is, since $\nabla(a)\in TS(\N)$, we can consider the fibre $(TZ)_{\nabla(a)}\subseteq TX(\N)$ of $d\rho$ over $\nabla(a)$.
We define this to be the {\em prolongation space of $Y=Z_a$}, and denote it by $\tau Y$.
That is, $\tau Y:=(TZ)_{\nabla(a)}$.

\begin{lemma}
The above definition of $\tau Y$ depends only on $Y$ and not on the presentation of $Y$ as $Z_a$.
\end{lemma}

\begin{proof}
Suppose $Y$ also appears as $Z'_b$ for some $0$-definable Zariski closed $Z'\subseteq S'\times X$ with $b\in S'(A)$ generic.
Replacing $b$ with $(a,b)$, we may assume that there is a dominant definable meromorphic map $f:S'\to S$ with $f(b)=a$, and that $Z'\subseteq Z\times_SS'$.
Hence $f':=(f,\id_X)\upharpoonright_Z':Z'\to Z$ restricts to the identity on $Z'_b=Y=Z_a$.
Moreover,  we have
$$\xymatrix{
Z\ar[d]_{\rho}& Z'\ar[l]_{f'}\ar[d]^{\rho'}\\
S& S'\ar[l]_{f}
}$$
which yields
$$\xymatrix{
TZ\ar[d]_{d\rho}& TZ'\ar[l]_{df'}\ar[d]^{d\rho'}\\
TS& TS'\ar[l]_{df}
}$$
Since $(A,\nabla)$ is a differential $\ccm$-structure, $df(\nabla b)=\nabla(a)$, so that
$$df'_{\nabla(b)}:(TZ')_{\nabla(b)}\to(TZ)_{\nabla(a)}$$
Since $df'=(df,\id_{TX})\upharpoonright_{T Z'}$, this shows that $(TZ')_{\nabla(b)}= (TZ)_{\nabla(a)}$, as desired.
\end{proof}

We denote the restriction of $\pi:TX\to X$ to $\tau Y$ also as $\pi:\tau Y\to Y$, and it is canonically attached to the prolongation space.
For any $b\in Y$, we denote the fibre by $\tau_bY$, and call it the {\em prolongation space to $Y$ at $b$}.
Note that if $Y$ is an $a$-definable Zariski closed subset of $X$ then $\tau Y$ is a $\nabla(a)$-definable Zariski closed subset of $TX$ and $\tau_b Y$ is a $\nabla(a)b$-definable Zariski closed subset of $T_bX$.

\begin{lemma}
\label{pronab}
Suppose $(B,\nabla)$ is a $\N$-valued differential $\ccm$-structure extending $(A,\nabla)$, and $b\in Y(B)$.
Then $\nabla(b)\in \tau_bY$.
\end{lemma}

\begin{proof}
Since $b\in X(B)$ we must have $\nabla(b)\in T_bX$.
So it remains to verify that $\nabla(b)\in\tau Y$.
Write $Y=Z_a$ as above.
Then $(\nabla(a),\nabla(b))=\nabla(a,b)\in T Z(\N)$.
In particular, $\nabla(b)\in (TZ)_{\nabla(a)}$ which is $\tau Y$ by construction.
\end{proof}

\begin{lemma}
\label{torsor}
If $Y$ is $A$-irreducible and $b\in Y$ is generic over $A$ then $\tau_bY$ is absolutely irreducible and $\dim(\tau_bY)=\dim Y$.
\end{lemma}

\begin{proof}
Let $a$ from $A$ be such that $Y=Z_a$ with $Z=\loc(a,b)\subseteq S\times X$ as above.
Because $\rho:Z\to S$ is dominant, $d\rho$ restricts to a surjective $\CC$-linear map between the tangent spaces at standard general points.
Hence, at the generic point in $\N$,  we have that $d\rho_{(a,b)}:T_{(a,b)}Z\to T_aS$ is a surjective $\CC(\N)$-linear map, where $\CC(\N)$ is the interpretation in $\N$ of the complex field, itself an algebraically closed field extending~$\CC$.
By definition, the tangent space $T_bY$ is the kernel of $d\rho_{(a,b)}$ while the prolongation space $\tau_bY$ is $d\rho_{(a,b)}^{-1}(\nabla a)$.
So $\tau_bY$ is a coset of $T_bY$ in $T_{(a,b)}Z$.
Absolute irreducibility of $\tau_bY$ follows, and $\dim(\tau_bY)=\dim(T_bY)$.
Finally, note that $\dim(T_bY)=\dim Y$ because for standard general $(p,q)\in Z(\M)$, the tangent space to $Z_p$ at $q$ is of dimension $\dim(Z_p)$.
\end{proof}

Finally, it is worth thinking about the case when $Y$ is $0$-definable.
That is, using the above notation, when $a\in\M$.
In that case, by Lemma~\ref{zerom}, $\nabla$ agrees with the zero section at $a$, and hence $\tau Y=TZ$ is just the tangent space of $Y$.
That is, for $0$-definable Zariski closed sets, the prolongation and tangent spaces agrees.

\bigskip
\section{Differentially closed $\ccm$-structures}
\label{sec:dccm}

\noindent
We aim to prove that $\ccmforallnabla$ admits a model companion.
We begin by exploring some properties of the existentially closed (e.c.) models.
This amounts to proving extension lemmas.
For example, Proposition~\ref{extdcl}, which says that every model of $\ccmforallnabla$ extends to the definable closure of the underlying model of $\ccm_\forall$, implies that if $(A,\nabla)$ is an e.c. model of $\ccmforallnabla$ then it is a differential $\ccm$-structure.
Moreover, the e.c. models of $\ccmforallnabla$ are precisely the existentially closed differential $\ccm$-structures.
This justifies:

\begin{definition}
A {\em differentially closed $\ccm$-structure} is an e.c. model of $\ccmforallnabla$.
\end{definition}

Here is the main extension lemma:

\begin{proposition}
\label{extccm}
Suppose $\N\models\ccm$ and $(A,\nabla)$ is a $\N$-valued differential $\ccm$-structure.
Suppose $X$ is an irreducible meromorphic variety, $b\in X(\N)$, and $Y:=\loc(b/A)$ is the smallest $A$-definable Zariski closed subset of $X(\N)$.
For any $c\in\tau_bY$ there is an extension of $\nabla$ to a $\N$-valued differential $\ccm$-structure on $\dcl(Ab)$ such that $\nabla(b)=c$.
\end{proposition}

\begin{proof}
Let $D:=\dcl(Ab)$.
We follow the approach of Proposition~\ref{extdcl}.
That is, given an element of $D$, say $d=f(a,b)$ where $a$ is from $A$ and $f:\loc(a,b)\to \loc(d)$ is a definable meromorphic map, we set  $\nabla(d):=df(\nabla a, c)$.
We have to verify that  $(\nabla a ,c)\in T_{(a,b)}\loc(a,b)$ for this to even make sense, that is to be able to apply $df$ to $(\nabla a,c)$.
Note, first of all, that since $Y=\loc(b/A)\subseteq\loc(b)$ we do have that $c\in\tau _bY\subseteq T_b\loc(b)$.
So $(\nabla a ,c)\in T_a\loc(a)\times T_b\loc(b)$, that is, $(\nabla a ,c)$ lies above $(a,b)$, and it only remains to check that $(\nabla a,c)\in T\loc(a,b)$.
Since  $Y=\loc(b/A)\subseteq \loc(b/a)$, and the latter is the fibre of the co-ordinate projection $\loc(a,b)\to\loc(a)$ over $a$, we have that $\tau Y\subseteq \tau\loc(b/a)$, and the latter is by definition the fibre of $T\loc(a,b)\to T\loc(a)$ over $\nabla(a)$.
Since $c\in\tau Y$, this tells us that $(\nabla a,c)\in T\loc(ab)$, as desired.

Considering the case when $d=a$ and $f:\loc(a,b)\to\loc(b)$ is the co-ordinate projection, we see that this definition of $\nabla$ on $D$ extends the given $\nabla$ on $A$.
Considering the case when $d=b$ (so that $a$ is the empty tuple and $f=\id$), we see that $\nabla(b)=c$, as desired.

While the proof of Proposition~\ref{extdcl} was carried out in the context of models of $\ccmforallnabla$, it works equally well in the setting of $\N$-valued differential $\ccm$-structures, showing that the way we have defined $\nabla$ on $D$ above yields, for any sort~$S$,  a well-defined map $\nabla:S(D)\to TS(\N)$ that is a section to $TS\to S$, and such that $dg(\nabla d)=\nabla(g(d))$ for any definable meromorphic map $g$ and tuple~$d\in S(D)$.
So $(D,\nabla)$ is again a $\N$-valued differential $\ccm$-structure.
\end{proof}

\begin{corollary}
\label{dccm-ccm}
If $(A,\nabla)$ is a differentially closed $\ccm$-structure then $A\models \ccm$.
\end{corollary}

\begin{proof}
We have that $A\subseteq\N$ for some $\N\models\ccm$.
Let $(B,\nabla)$ be a maximal $\N$-valued differential $\ccm$-structure extending $(A,\nabla)$.
This exists as $\N$-valued differential $\ccm$-structures are preserved under unions of chains, as can be easily verified from the definition.

We claim that $B=\N$.
Given $b\in X(\N)$ for some sort $X$, let $Y=\loc(b/B)$ and choose $c\in\tau_bY$.
By Proposition~\ref{extccm} we can extend $\nabla$ to a $\N$-valued differential $\ccm$-structure on $\dcl(Bb)$.
By maximality, it follows that $b\in X(B)$ to start with.
As $X$ and $b$ were arbitrary, this shows that $B=\N$.

We have that $(A,\nabla)\subseteq(\N,\nabla)$ is an extension of differential $\ccm$-structures.
By quantifier elimination, $\ccm$ has a universal-existential axiomatisation.
Since $(A,\nabla)$ is existentially closed, the truth of such axioms in $(\N,\nabla)$ will imply their truth in $(A,\nabla)$.
That is, $A\models\ccm$, as desired.
\end{proof}

This is, of course, not enough.
That is, not every differential $\ccm$-structure on a model of $\ccm$ is differentially closed.
For example, the standard model $\M$ admits the trivial differential structure $\nabla=0$, but is not existentially closed as we can use Proposition~\ref{extccm} to produce nontrivial differential $\ccm$-structures extensions.

The following property of differentially closed $\ccm$-structures, which we refer to as the {\em geometric axiom}, can be read as saying that $\nabla$ is a ``generic" section to the tangent space:

\begin{proposition}
\label{ga}
If $(\N,\nabla)$ is a differentially closed $\ccm$-structure then it satisfies the following condition:
\begin{itemize}
\item[($\ga$)]
Suppose $S$ is a sort, $X\subseteq S$ is an $\N$-definable irreducible Zariski closed subset, $Y\subseteq \tau X$ is an $\N$-definable irreducible Zariski closed subset that projects dominantly onto $X$, and $Y_0\subsetneq Y$ is a proper $\N$-definable Zariski closed subset.
Then there exists $a\in X(\N)$ such that $\nabla(a)\in Y\setminus Y_0$.
\end{itemize}
\end{proposition}

\begin{proof}
We already know, by Corollary~\ref{dccm-ccm}, that $\N\models\ccm$.
Let $\U\succeq\N$ be a sufficiently saturated elementary extension, and let $c\in Y(\U)$ be generic in $Y$ over~$\N$.
In particular, $c\in Y\setminus Y_0$.
By dominance, $b:=\pi(c)\in X(\U)$ is generic over $\N$.
In particular, $\loc(b/\N)=X$ and $c\in\tau_bX(\U)$.
So, by Proposition~\ref{extccm}, we can extend $\nabla$ to a $\U$-valued differential $\ccm$-structure on $\dcl(\N b)$ such that $\nabla(b)=c$.
Then, as in the proof of Corollary~\ref{dccm-ccm}, we can extend $\nabla$ further to all of $\U$ so that $(\U,\nabla)\models\ccmforallnabla$.
Now, $b$ witnesses that in $(\U,\nabla)$ there is a point of $X$ that is sent by $\nabla$ into $Y\setminus Y_0$.
By existential closedness of $(\N,\nabla)$, there must exist $a\in X(\N)$ such that $\nabla(a)\in Y\setminus Y_0$.
\end{proof}

As the terminology already indicates, the geometric axiom characterises differentially closed $\ccm$-structures:

\begin{theorem}
\label{dccm}
A model $(\N,\nabla)\models\ccmforallnabla$ is existentially closed if and only if $\N\models\ccm$ and condition~$(\ga)$ of Proposition~\ref{ga} holds.
\end{theorem}

\begin{proof}
Corollary~\ref{dccm-ccm} and Proposition~\ref{ga} gave the left-to-right direction.
We therefore assume that $\N\models\ccm$ and $(\N,\nabla)$ satisfies~$(\ga)$, and show that $(\N,\nabla)$ is existentially closed.
Let $S$ be a sort, $x$ a variable belonging to $S$, and $\phi(x)$ a (finite) conjunction of $L_\nabla$-literals over $\N$ that is realised by $c\in S(A)$ in some extension $(A,\nabla)\models\ccmforallnabla$ of $(\N,\nabla)$.
We need to show that $\phi(x)$ has a realisation already in $(\N,\nabla)$.
As in the proof of Corollary~\ref{dccm-ccm}, we can extend $(A,\nabla)$ further to $(\U,\nabla)\models\ccmforallnabla$ where $\U\models\ccm$.

Let $d$ be the {\em order} of $\phi(x)$, that is, the largest positive integer such that $\nabla^d(x)$, namely $\nabla$ iterated $d$-times and applied to $x$, appears in $\phi(x)$.
We leave it to the reader to verify that $\phi(x)$ can then be rewritten as $\big(\nabla^d(x)\in U\big)\wedge\big(\nabla^d(x)\notin V\big)$ where $U$ and $V$ are $\N$-definable Zariski closed subsets of $T^d(S)$, the $d$th iterated tangent space of $S$.

Let 
$Y:=\loc(\nabla^d c/\N)\subseteq T^d(S)(\U)$.
Since $c$ realises $\phi(x)$, we must have that $\nabla^d(c)\in U\setminus V$ and so $Y\subseteq U\setminus V$.
In particular, $Y_0:=Y\cap V$ is a proper $\N$-definable Zariski closed subset of $Y$.
We aim to find $a\in S(\N)$ such that $\nabla^d(a)\in Y\setminus Y_0$; this will suffice as such an $a$ would be a realisation of $\phi(x)$ in $(\N,\nabla)$.

Let $\overline c:=\nabla^{d-1}(c)$ and $X:=\loc(\overline c/\N)$.
Then $\nabla^d(c)=\nabla(\overline c)$, so that $Y$ is contained in $\tau X$ and projects dominantly onto $X$.
Hence, by~$(\ga)$, there is an $\overline a\in X(\N)$ such that $\nabla(\overline a)\in Y\setminus Y_0$.
Consider the first co-ordinate projection $\pi:T^{d-1}(S)\to S$, and set $a:=\pi(\overline a)\in S(\N)$.
It will suffice to show, therefore, that $\nabla^{d-1}(a)=\overline a$.

For each $\ell\geq 0$, let us denote by $\pi_\ell: T^{\ell+1}S\to T^\ell S$ the canonical projection.
Moreover, for each $\ell=0,\dots,d-1$, let $\overline a_\ell$ be the image of $\overline a$ in $T^\ell S$.
So, in particular, $\overline a_0=a$ and $\overline a_{d-1}=\overline a$.
We claim that it suffices to show that
\begin{equation}
\label{a}
\overline a_{\ell+1}=\nabla(\overline a_\ell)
\end{equation}
for all $\ell=0,\dots,d-2$.
Indeed, this would imply that
$\overline a =\overline a_{d-1}=\nabla(\overline a_{d-2})=\nabla^2(\overline a_{d-3})=\cdots=\nabla^{d-1}(a)$,
as desired.
So let us fix $\ell=0,\dots,d-2$ and show~$(\ref{a})$.
The idea is to construe~$(\ref{a})$ as a Zariski closed condition on $\nabla(\overline a)$.
First of all, noting that 
$\pi_{\ell+1}(\nabla \overline a_{\ell+1})=\overline a_{\ell+1}$
and
$d\pi_\ell(\nabla a_{\ell+1})=\nabla( \pi_\ell a_{\ell+1})=\nabla(\overline a_\ell)$,
we see that~$(\ref{a})$ is equivalent to
\begin{equation}
\label{b}
\pi_{\ell+1}(\nabla \overline a_{\ell+1})=d\pi_\ell(\nabla \overline a_{\ell+1}).
\end{equation}
Next, letting $\rho:T^{d-1}S\to T^{\ell+1}S$ be the projection, we have that $\rho(\overline a)=\overline a_{\ell+1}$, and hence
$\nabla(\overline a_{\ell+1})=\nabla(\rho \overline a)=d\rho(\nabla(\overline a))$.
So~$(\ref{b})$ is equivalent to
\begin{equation}
\label{c}
\pi_{\ell+1}d\rho(\nabla \overline a)=d(\pi_\ell\rho)(\nabla \overline a).
\end{equation}
This is a Zariski closed condition on $\nabla\overline a$, and as $\nabla \overline a$ is in $Y=\loc(\nabla\overline c/\N)$, it suffices to verify that the identity holds of $\nabla\overline c$.
But this follows from the fact that $\nabla \overline c=\nabla^d c$,
\begin{eqnarray*}
\pi_{\ell+1}d\rho(\nabla \overline c)
&=&
\pi_{\ell+1}d\rho(\nabla^d c)\\
&=&
\pi_{\ell+1}\nabla(\rho(\nabla^{d-1} c))\\
&=&
\rho(\nabla^{d-1} c)\\
&=&
\nabla^{\ell+1}c\\
&=&
\nabla(\nabla^{\ell} c)\\
&=&
\nabla(\pi_\ell\rho(\nabla^{d-1}c))\\
&=&
d(\pi_\ell\rho)(\nabla^dc))\\
&=&
d(\pi_\ell\rho)(\nabla\overline c)).
\end{eqnarray*}
Hence~$(\ref{c})$ holds, as desired.
\end{proof}

That condition~$(\ga)$ of Proposition~\ref{ga} is first-order expressible follows from the fact that as $X$ varies in an $L$-definable family, $\tau X$ varies in an $L_{\nabla}$-definable family by construction (see Section~\ref{prolongations}), and that in $\ccm$ irreducibility and domination are definable in parameters (see~\cite[Section~2]{ret}).
Theorem~\ref{dccm} thus gives us a model companion to $\ccmforallnabla$, namely the theory of differentially closed $\ccm$-structures, which we denote $\dccm$.

\bigskip
\section{Basic model theory of $\dccm$}
\label{bmt}

\noindent
From general model theory, we have that $\dccm$ is model-complete.
In this section we prove that $\ccmforallnabla$ has the amalgamation property, from which we can deduce that $\dccm$ is complete and admits quantifier elimination.
As a consequence we obtain a geometric description of algebraic and definable closure.

But first we need an extension lemma for algebraic closure, whereas we have only dealt with definable closure (in Proposition~\ref{extdcl}) so far.

\begin{lemma}
\label{extacl}
Suppose $(A,\nabla)$ is a differential $\ccm$-structure with $A\subseteq\N\models\ccm$, and $b\in\acl(A)$.
Then there is a unique $\N$-valued differential $\ccm$-structure on $\dcl(Ab)$ that extending $\nabla$.
Moreover, this extension is in fact $\dcl(Ab)$-valued.
\end{lemma}

\begin{proof}
Let $Y=\loc(b/A)$ and $c\in\tau_bY$.
By Proposition~\ref{extccm} we can extend $\nabla$ from $A$ to a $\N$-valued differential $\ccm$-structure on  $\dcl(Ab)$ by sending $\nabla(b):=c$.
We will show that $\tau_bY=\{c\}$ is a singleton and hence $c\in\dcl(Ab)$, so that the above extension is in fact $\dcl(Ab)$-valued, and so $(\dcl(Ab),\nabla)\models\ccmforallnabla$.
This will also show uniqueness as any extension of $\nabla$ to $\dcl(Ab)$ would have to take $b$ into $\tau_bY=\{c\}$, by Lemma~\ref{pronab}, and hence would agree with the one we just constructed.

Let $X:=\loc(b)$ and write $Y=Z_a$ where $Z=\loc(a,b)\subseteq S\times X$ is a $0$-definable irreducible Zariski closed set and $S$ is a sort with $a\in S(A)$ generic.
The fact that $b\in\acl(A)$ means that $Y$ is finite, and hence the co-ordinate projection $\rho:Z\to S$ is generically finite-to-one.
It follows that $d_p\rho:T_pZ\to T_{\rho(p)}S$ is an isomorphism for general $p\in Z(\M)$.
Hence $d_{(a,b)}\rho:T_{(a,b)}Z\to T_aS$ is a bijection.
If  $c,c'\in\tau_bY$ then we know, by the proof of Proposition~\ref{extccm},  that $(\nabla a,c), (\nabla a,c')\in T_{(a,b)}Z$.
But as $d\rho$ takes both $(\nabla a,c)$ and $(\nabla a,c')$ to $\nabla(a)\in T_aS$ we must have $c=c'$.
So $\tau_bY=\{c\}$, as desired.
\end{proof}

Next we prove independent amalgamation.
We will use $\displaystyle \ind^{\ccm}$ to mean nonforking independence in $\ccm$.

\begin{lemma}
\label{indamalgam}
Suppose $A,B_1,B_2$ are definably closed subsets of $\N\models \ccm$, with $A\subseteq B_1\cap B_2$ and $\displaystyle B_1\ind_A^{\ccm} B_2$.
Suppose $\nabla_i$ is a differential $\ccm$-structure on $B_i$, for $i=1,2$, such that $\nabla_1$ and $\nabla_2$ agree on $A$.
Then there is a common extension, $\nabla$, of $\nabla_1$ and $\nabla_2$ to $B:=\dcl(B_1B_2)$ such that $(B,\nabla)\models\ccmforallnabla$.
\end{lemma}

\begin{proof}
Using Lemma~\ref{extacl} we can extend the differential $\ccm$-structure on $A,B_1,B_2$ uniquely to their algebraic closures in $\N$.
In particular, $\nabla_1$ and $\nabla_2$ will agree on $\acl(A)$.
So we may as well assume that $A=\acl(A)$, and $B_i=\acl(B_i)$ for $i=1,2$.
One consequence of $A$ being $\acl$-closed is that Zariski loci over $A$ are absolutely irreducible, and hence independence over $A$ has the following Zariski-topological characterisation:
$$\displaystyle b_1\ind_A^{\ccm}b_2\ \text{ if and only if }\ \loc(b_1,b_2/A)=\loc(b_1/A)\times\loc(b_2/A).$$
See~\cite[Section~2]{ret} for details.

Every tuple from $B$ is of the form $b=f(b_1,b_2)$ where each $b_i$ is from $B_i$, and $f:\loc(b_1,b_2)\to\loc(b)$ is a definable meromorphic map.
Our only choice is to define $\nabla(b):=df(\nabla_1 b_1,\nabla_2 b_2)$.
But we need $(\nabla_1 b_1,\nabla_2 b_2)\in T_{(b_1,b_2)}\loc(b_1,b_2)$ for this to make sense.
This is what we now check.

Let $a$ be a tuple from $A$ such that $\loc(b_1,b_2/A)=\loc(b_1,b_2/a)$.
Let us denote by $\nabla$ the common restriction of $\nabla_1$ and $\nabla_2$ to $A$.
Taking prolongations with respect to the differential $\ccm$-structure $(A,\nabla)$, and using that for $i=1,2$ we have $(A,\nabla)\subseteq(B_i,\nabla_i)$, we see that $\nabla_i(b_i)\in\tau_{b_i}\loc(b_i/a)$.
Hence,
\begin{eqnarray*}
(\nabla_1 b_1,\nabla_2 b_2)
&\in&
\tau_{b_1}\loc(b_1/a)\times \tau_{b_2}\loc(b_2/a)\\
&=&
\tau_{(b_1,b_2)}\big(\loc(b_1/a)\times\loc(b_2/a)\big)\\
&=&\tau_{(b_1,b_2)}\loc(b_1,b_2/a)\ \ \ \ \ \ \ \ \ \text{ as $b_1\ind^{\ccm}_a b_2$}\\
&\subseteq &
T_{(b_1,b_2)}\loc(b_1,b_2)
\end{eqnarray*}
as desired.

So it does make sense to set $\nabla(b):=df(\nabla_1 b_1,\nabla_2 b_2)$ for $b=f(b_1,b_2)$.
The next step is to make sure this is well defined.
What if we also have $b=f'(b_1',b_2')$?
This is dealt with exactly as it was done in Proposition~\ref{extdcl}.
Namely, let $Z:=\loc(b_1,b_2,b_1',b_2')$ and consider $\bar f:=(f,f'):Z\to \loc(b)^2$.
Since $\bar f$ takes a generic point of $Z$ to the diagonal we have that $d\bar f:TZ\to T(\loc(b)^2)=(T\loc(b))^2$ lands in the diagonal.
Now, the argument in the previous paragraph, applied to $b_ib_i'$, shows, in particular, that $(\nabla_1 (b_1b_1'),\nabla_2 (b_2b_2'))\in T\loc(b_1b_1',b_2b_2')$.
Hence, $(\nabla_1b_1,\nabla_2b_2,\nabla_1b_1',\nabla_2b_2')\in TZ$ and we get that $df(\nabla_1 b_1,\nabla_2 b_2)=df'(\nabla_1 b_1',\nabla_2 b_2')$.

We have defined $\nabla$ on $B$ in such a way that it is a section to $TS\to S$ for any sort~$S$.
It remains to check Axiom~$(2)$ of Definition~\ref{ccmforallnabla}.
That is, suppose $g:S\to S'$ a definable meromorphic map between sorts, and $b\in S(B)$, $b'\in S'(B)$ with $g(b)=b'$.
We need to show that $dg(\nabla b)=\nabla b'$.
We may assume that there are $b_1, b_2$ from $B_1,B_2$, respectively, and definable meromorphic maps $f, f'$ such that 
$b=f(b_1,b_2)$ and $b'=f'(b_1,b_2)$.
It follows that $gf=f'$ on $\loc(b_1,b_2)$, and so we compute:
\begin{eqnarray*}
dg(\nabla b)
&=& dg(df(\nabla_1 b_1,\nabla_2 b_2))\ \ \text{ by definition of $\nabla(b)$}\\
&=& d(gf)(\nabla_1 b_1,\nabla_2 b_2))\ \ \text{ by functoriality}\\
&=&df'(\nabla_1 b_1,\nabla_2 b_2))\ \ \ \text{ as $gf=f'$}\\
&=& \nabla b'\ \ \ \ \text{ by definition of $\nabla(b')$.}
\end{eqnarray*}
This completes the proof that $(B,\nabla)\models\ccmforallnabla$.
\end{proof}

\begin{proposition}
\label{apqe}
$\ccmforallnabla$ has the amalgamation property.
In particular, $\dccm$ admits quantifier elimination and is complete.
\end{proposition}

\begin{proof}
Suppose $(B_i,\nabla)\models\ccmforallnabla$, for $i=1,2$, with a common substructure $(A,\nabla)$.
We seek a model $(B,\nabla)\models\ccmforallnabla$ into which $(B_1,\nabla)$ and $(B_2,\nabla)$ both embed over $A$.
Let $\U\supseteq B_1$ be a sufficiently saturated model of $\ccm$.
By universality, there is an embedding of $B_2$ into $\U$ over $A$.
Moreover, after taking nonforking extensions in $\ccm$, we can find such an embedding such that the image of $B_2$ is independent from $B_1$ over $A$.
We may as well assume, therefore, that $\displaystyle B_2\subseteq\U$ already, and that $\displaystyle B_1\ind_A^{\ccm}B_2$.
Applying Lemma~\ref{indamalgam}, we have a differential $\ccm$-structure $\nabla$ on $B:=\dcl(B_1B_2)$ that extends $\nabla$ on both both $B_1$ and $B_2$.

Quantifier elimination now follows for $\dccm$, as a general consequence for a model companion of a universal theory with amalgamation.

Completeness also follows as we have a prime substructure: all differentially closed $\ccm$-structures extend the standard model $\M\models\ccm$ equipped with the trivial differential structure (Lemma~\ref{zerom}).
\end{proof}

Next, we wish to characterise definable and algebraic closure in $\dccm$.
First of all, given $(\N,\nabla)\models\dccm$ and $A\subseteq \N$, let us denote by $\langle A\rangle$ the $L_\nabla$-structure generated by $A$.
If $A$ is already an $L_\nabla$-substructure and $a$ is a tuple then we denote by $A\langle a\rangle$ the $L_\nabla$-structure generated by $A\cup\{a\}$.
Note that
$$A\langle a\rangle=A\cup\{a,\nabla(a), \nabla^2(a),\dots\}.$$
Quantifier elimination tells us that $\tp(a/A)=\tp(a'/A)$ if and only if there is an $L$-isomorphism $\alpha:A\langle a\rangle\to A\langle a'\rangle$ that fixes $A$ point-wise and sends $\nabla^n(a)$ to $\nabla^n(a')$ for all $n\geq0$.

We have been using $\acl$ and $\dcl$ for algebraic and definable closure in the $L$-theory $\ccm$.
We will continue to do so, using $\aclnabla$ and $\dclnabla$ for algebraic and definable closure in the $L_\nabla$-theory $\dccm$.

\begin{proposition}
\label{dclacl}
Suppose $(\N,\nabla)$ is a differentially closed $\ccm$-structure and $A\subseteq \N$.
Then $\dclnabla(A)=\dcl(\langle A\rangle)$ and $\aclnabla(A)=\acl(\langle A\rangle)$.
\end{proposition}

\begin{proof}
By Proposition~\ref{extdcl}, $\dcl(\langle A\rangle)$ is a differential $\ccm$-substructure of $(\N,\nabla)$.
Replacing $A$ by $\dcl(\langle A\rangle)$, we may as well assume that $A$ is a differential $\ccm$-substructure to start with, and show that
$\dclnabla(A)= A$ and  $\aclnabla(A)= \acl(A)$.
The right-to-left containments are clear.

For the converses, let us first suppose that $b\notin \acl(A)=:B$.
By Lemma~\ref{extacl}, $(B,\nabla)$ is a differential $\ccm$-substructure of $(\N,\nabla)$.
We can find, in some $\N\succeq\N$, a copy of $\N$ over $B$, say $\N'$, such that $\displaystyle \N\ind^{\ccm}_B\N'$.
Let $\alpha:\N\to \N'$ be an $L$-isomorphism over $B$ witnessing this, and consider $b':=\alpha(b')$.
The fact that $\displaystyle b\ind^{\ccm}_Bb'$ and that $b\notin B$ forces $b\neq b'$.
On the other hand, setting $\nabla':=\alpha\nabla\alpha^{-1}$ we have that $(\N',\nabla')\models\dccm$ and that $\alpha:(\N,\nabla)\to (\N',\nabla')$ is an $L_{\nabla}$-isomorphism over $B$.
Now, we can find a common extension of $\nabla$ and $\nabla'$ to $\dcl(\N\N')$ in $\N$ by Lemma~\ref{indamalgam} and then further to a model $(\K,\nabla)\models\dccm$.
So, in $(\K,\nabla)$ we have produced at least two distinct realisations, $b$ and $b'$, of $\tp(b/B)$.
Repeating the process we can show that $\tp(b/B)$ has arbitrarily many realisations.
That is, $b\notin\aclnabla(B)=\aclnabla(A)$, as desired.

Finally, suppose, toward a contradiction, that  $b\in\dclnabla(A)\setminus A$.
This time we produce two distinct realisations of $\tp(b/A)$ for our contradiction.
Since $\dclnabla(A)\subseteq\aclnabla(A)=\acl(A)$, we have that $b\in\acl(A)\setminus A$.
Hence $\tp_L(b/A)$ has a second realisation, $b'\in\acl(A)$ with $b'\neq b$.
We thus have an $L$-isomorphism $\alpha:\dcl(Ab)\to\dcl(Ab')$ that fixes $A$ point-wise and sends $b$ to $b'$.
But, by Lemma~\ref{extacl}, $\dcl(Ab)$ and $\dcl(Ab')$ are differential $\ccm$-substructures of $(\N,\nabla)$, and, as they each admit unique differential structures extending $\nabla$ on $A$, we must have that $\alpha$ is an $L_{\nabla}$-isomorphism.
By quantifier elimination, this means $\tp(b/A)=\tp(b'/A)$.
\end{proof}

\bigskip
\section{Stability and elimination of imaginaries}
\label{sec:stability}

\noindent
We work now in a fixed sufficiently saturated model $(\U,\nabla)\models\dccm$ and adopt the usual convention that all parameter sets are assumed to be of cardinality less than that of the saturation.

In order to prove that $\dccm$ is a stable theory, and to capture the meaning of nonforking independence therein, we will follow an axiomatic approach.
That is, we first introduce a natural notion of independence and then show that it has all the properties that characterise nonforking independence in stable theories.

\begin{definition}
Given sets $A,B,C$, we  say that {\em $A$ is  independent of $B$ over $C$}, denoted by
$\displaystyle A\ind_CB$,
to mean that
$\displaystyle \dclnabla(A)\ind^{\ccm}_{\dclnabla(C)}\dclnabla(B)$.
\end{definition}

Note that we do not yet know that $\ind$ is nonforking independence, but we allow ourselves the notation as we will soon see that it is.

Let us first verify that $\ind$ is a notion of independence, in the sense introduced by Kim and Pillay~\cite{kimpillay}.
First of all, it is clearly invariant under the action of automorphisms of $(\U,\nabla)$.
Local character, finite character, symmetry, and transitivity all follow easily from the corresponding properties for  $\displaystyle \ind^{\ccm}$.

\begin{lemma}[Extension]
Given $a, C\subseteq B$ there is $a'\models\tp(a/C)$ such that $\displaystyle a'\ind_CB$.
\end{lemma}

\begin{proof}
We may assume that $C=\dclnabla(C)$ and $B=\dclnabla(B)$.
By extension in $\ccm$ there is sequence $\displaystyle (a_n':n\geq0)\ind^{\ccm}_CB$ and an $L$-isomorphism
$$\alpha:C\langle a\rangle\to C\cup\{a_n':n\geq0\}$$
that fixes $C$ point-wise and takes $\nabla^n(a)$ to $a'_n$ for all $n\geq0$.
Extend $\alpha$ to an $L$-isomorphism
$$\alpha: A:=\dcl(C\langle a\rangle)\to A':=\dcl(C\cup\{a_n':n\geq0\}).$$
Set $\nabla':=\alpha\nabla\alpha^{-1}$ on $A'$ so that $(A',\nabla')$ is a differential $\ccm$-structure isomorphic to $(A,\nabla)$.
On the other hand, since $\displaystyle A'\ind^{\ccm}_CB$, Lemma~\ref{indamalgam} gives us a common extension of $(A',\nabla')$ and $(B,\nabla)$ to a model $(\N,\nabla')\models \dccm$.
By universality we have an embedding  $\iota:(\N,\nabla')\to (\U,\nabla)$ over $B$.
Then
$\beta:=\iota\circ\alpha: A \to\U$
is an $L$-isomorphism with its image that fixes $C$ point-wise and takes $\nabla^n(a)$ to $\nabla^n(\beta(a))$ for all $n\geq0$.
Hence, by quantifier elimination, $a':=\beta(a)\models\tp(a/C)$.
Now $\displaystyle A'\ind^{\ccm}_CB$ implies that $\displaystyle \iota(A')\ind^{\ccm}_CB$ as $\iota$ is over $B$.
But $\iota(A')=\beta(A)=\dcl(C\langle a'\rangle)$.
So  $\displaystyle a'\ind_CB$, as desired.
\end{proof}

\begin{lemma}[Stationarity over algebraically closed sets]
Suppose $C=\aclnabla(C)\subseteq B$, and $a,a'$ are tuples.
If $\tp(a/C)=\tp(a'/C)$, and both $a$ and $a'$ are independent of $B$ over $C$, then $\tp(a/B)=\tp(a'/B)$.
\end{lemma}

\begin{proof}
We may assume, without loss of generality, that $B=\dclnabla(B)$.
Since $\tp(a/C)=\tp(a'/C)$, there is an $L$-isomorphism $\alpha:C\langle a\rangle\to C\langle a'\rangle$ that fixes $C$ pointwise and takes $\nabla^n(a)$ to $\nabla^n(a')$ for all $n\geq 0$.
Since $C$ is $\acl$-closed, and the sequences $(\nabla^n a:n\geq0)$ and $(\nabla^n a':n\geq0)$ are both $\ccm$-independent from $B$ over~$C$, stationarity over algebraically closed sets in $\ccm$ implies that there is an $L$-isomorphism $\beta: B\langle a\rangle\to B\langle a'\rangle$ that fixes $B$ pointwise and takes $\nabla^n(a)$ to $\nabla^n(a')$ for all $n\geq 0$.
By quanitifier elimination, $\tp(a/B)=\tp(a'/B)$.
\end{proof}

\begin{corollary}
\label{stable}
$\dccm$ is a stable theory and $\ind$ is nonforking independence.
\end{corollary}

\begin{proof}
This follows from the above observations by the characterisation of nonforking independence in simple (and hence stable) theories, see~\cite{kimpillay}.
\end{proof}

We can also deduce stability by counting types. 
In fact, we get total transcendentality:

\begin{theorem}
\label{lambdastable}
$\dccm$ is $\lambda$-stable for every cardinal $\lambda\geq 2^{\aleph_0}$.
\end{theorem}

\begin{proof}
We count types.
Fix $\lambda\geq 2^{\aleph_0}$ and a subset $A\subseteq\U$ of cardinality at most $\lambda$.
We show that there are at most $\lambda$-many complete types over $A$.
We may assume that $A=\dclnabla(A)$ is a differential $\ccm$-substructure.

Suppose $X$ is a sort,  $a\in X(\U)$, and consider $\tp(a/A)$.
By quantifier elimination, it is determined by the sequence of types $\big(\tp_L(\nabla^n a/A):n\geq0\big)$ in $\ccm$.
Let $$Z_n:=\loc(\nabla^n a/A\nabla^{n-1}a).$$
I claim that there is some $N\geq 0$ such that
$Z_{n+1}=\tau_{\nabla^na}(Z_n)$, for all $n\geq N$.
This will suffice, since then $\tp(a/A)$ is determined by the pair $\big(N,\tp_L(\nabla^N a/A)\big)$ of which there are at most $\lambda$-many possibilities by the $\lambda$-stability of $\ccm$.

Note that $\nabla^{n+1}(a)\in\tau_{\nabla^n a}(Z_n)$, and so  $Z_{n+1}\subseteq \tau_{\nabla^n a}(Z_n)$.
But $\dim(\tau_{\nabla^n a}(Z_n))=\dim(Z_n)$ by Lemma~\ref{torsor}.
So $\dim(Z_{n+1})$ is a nonincreasing function of $n$ that must eventually stabilise.
By irreducibility of $\tau_{\nabla^n a}(Z_n)$, this forces $Z_{n+1}= \tau_{\nabla^n a}(Z_n)$ for large enough~$n$, as desired.
\end{proof}

\begin{theorem}
\label{ei}
$\dccm$ admits elimination of imaginaries.
\end{theorem}

\begin{proof}
A general criterion for elimination of imaginaries in a stable theory is that finite sets have codes and global types have canonical bases, in the home sorts; see, for example,~\cite[Section~3]{johnson}.
That finite sets in $\dccm$ have codes in the home sort follows from elimination of imaginaries in $\ccm$, see~\cite{pillay00} and~\cite[Appendix]{saturated}.

So, we fix a saturated $\N\preceq\U$ and a complete type $p=\tp(a/\N)$, and show that~$p$ has a canonical base in $\N$.
Let $N$ be as in the proof of Theorem~\ref{lambdastable}; that is,
$$\loc(\nabla^{n+1} a/\N\nabla^na)=\tau_{\nabla^na}\big(\loc(\nabla^n a/\N\nabla^{n-1}a)\big),$$
for all $n\geq N$.
By elimination of imaginaries for $\ccm$, there is a code $c$ for $\loc(\nabla^Na/\N)$ in $\N$.
We claim that $c$ is a canonical base for $p$.

Fix $\sigma$, an $L_\nabla$-automorphism of $\N$.
We need to show that $p^\sigma=p$ if and only if $\sigma(c)=c$.
One direction is clear: if $p^\sigma=p$ then
$\tp_L(\nabla^Na/\N)^\sigma=\tp_L(\nabla^Na/\N)$
and hence
$\loc(\nabla^Na/\N)^\sigma=\loc(\nabla^Na/\N)$,
so that $\sigma(c)=c$.

For the converse, suppose $\sigma(c)=c$.
Extend $\sigma$ to an $L_\nabla$-automorphism $\hat\sigma$ of $\U$, and let $\hat a:=\hat\sigma(a)$.
We have that
\begin{eqnarray*}
\tp_L(\nabla^N\hat a/\N)
&=& \tp_L(\nabla^N a/\N)\ \ \ \ \ \ \ \ \text{ since $\sigma(c)=c$}\\
\loc(\nabla^{N+1} a/\N\nabla^Na)
&=&\tau_{\nabla^Na}\big(\loc(\nabla^N a/\N\nabla^{N-1}a)\big)\ \text{ by choice of $N$, and}\\
\loc(\nabla^{N+1} \hat a/\N\nabla^N\hat a)
&=&\tau_{\nabla^N\hat a}\big(\loc(\nabla^N \hat a/\N\nabla^{N-1}\hat a)\big)\ \text{ by applying $\hat\sigma$.}
\end{eqnarray*}
These imply that $\tp_L(\nabla^{N+1}\hat a/\N)=\tp_L(\nabla^{N+1} a/\N)$.
We can iterate to prove that $\tp_L(\nabla^{n}\hat a/\N)=\tp_L(\nabla^{n} a/\N)$ for all $n\geq 0$.
By quantifier elimination,
$$p^\sigma=\tp(\hat a/\N)=\tp(a/\N)=p,$$
as desired.
\end{proof}

\bigskip
\section{Meromorphic vector fields and finite-dimensional types}
\label{fd}

\noindent
We return in this final section to the motivating objects of interest: meromorphic vector fields.
Our goal is to show that they are captured, up to bimeromorphic equivalence, in $\dccm$, by the ``finite-dimensional" types.

We continue to work in a fixed sufficiently saturated model $(\U,\nabla)\models\dccm$.

\begin{definition}
Suppose $A$ is an $L_\nabla$-substructure and $p=\tp(b/A)$ is a complete type.
By the {\em dimension of $p$}, denote by $\dimnabla(p)$ or $\dimnabla(b/A)$, we mean the sequence of nonnegative nondecreasing integers $\big(\dim(\loc(\nabla^nb/A)):n\geq 0\big)$ ordered lexicographically.
If $\dimnabla(p)$ is eventually constant then we say that $p$ is {\em finite-dimensional} and we (re)use $\dimnabla(p)$ to denote that eventual finite number.
\end{definition}

Note that the dimension depends only on the type $p$ and not on the choice of realisation $b$.
On the other hand, this dimension is not invariant under definable bijection -- for example $b$ and $\nabla(b)$ are interdefinable over the empty set but the dimension sequences are not always the same (one is a shift of the other).
Nevertheless, whether or not a type is finite-dimensional, and the value of that finite dimension in the case that it is, is invariant under definable bijection.

Dimension witnesses forking:

\begin{proposition}
\label{nfdim}
Suppose $a$ is a tuple and $C\subseteq B$ are $L_\nabla$-substructures.
Then $\displaystyle a\ind_CB$ if and only if $\dimnabla(a/B)=\dimnabla(a/C)$.
\end{proposition}

\begin{proof}
We may assume that $B=\dclnabla(B)$ and $C=\dclnabla(C)$.
We have, by Proposition~\ref{dclacl}, that $\dclnabla(Ca)=\dcl(C\langle a\rangle)$.
Hence $\displaystyle a\ind_CB$ is equivalent to $\displaystyle \nabla^n(a)\ind^{\ccm}_CB$ for all $n\geq 0$.
But, as dimension witnesses forking in $\ccm$, this is equivalent to $\dim(\loc(\nabla^n a/B))=\dim(\loc(\nabla^n a/C))$ for all $n\geq 0$.
\end{proof}

It follows that if $p$ is finite-dimensional then it is of finite $U$-rank, and in fact that $U(p)\leq\dimnabla(p)$.

A natural source of finite-dimensional types over the empty set are meromorphic vector fields in the sense of Definition~\ref{mervf}.
Suppose $(X,v)$ is such.
Consider the type $p(x)$, over the empty, which says that $x\in X$ is generic and that $\nabla(x)=v(x)$.
This is consistent by the geometric axiom of Proposition~\ref{ga}.
Indeed, given any proper Zariski closed $X_0\subseteq X$, apply~$(\ga)$ to $Y$ the Zariski closure of the image of $v$ in $TX=\tau X$ and $Y_0$ the restriction of $Y$ to $X_0$, yielding a $\U$-point $a\in X\setminus X_0$ with $\nabla(a)=v(a)$.
Moreover, by quantifier elimination, this type is complete: the $L$-type of $x$ is determined by $x$ being generic in $X$, and $\nabla(x)=v(x)$ implies $\nabla^n(x)=v_n(x)$ for appropriate  definable meromorphic $v_n:X\to T^nX$, for all $n$.
We call $p$ the {\em generic type of $(X,v)$}.
Note that $p$ is finite-dimensional;  in fact, $\dimnabla(p)=\dim X$.
Indeed, for all $n\geq 0$, we have that $\nabla^n(b)=v_n(b)$ and $v_n$ is a definable meromorphic section to $T^nX\to X$, and hence, as $b$ is generic in $X$, we get $\dim(\loc(\nabla^nb/A))=\dim X$.

It turns out that all finite-dimensional types arise this way:

\begin{theorem}
\label{fdmervf}
Every finite-dimensional type over the empty set in $\dccm$ is, up to interdefinability, the generic type of a meromorphic vector field.
\end{theorem}

\begin{proof}
Suppose $p=\tp(b)$ is finite dimensional.
Let $d\geq 0$ be such that
$$\dim(\loc(\nabla^{d+1}b))=\dim(\loc(\nabla^{d}b)).$$
Since the projection $\loc(\nabla^{d+1}b)\to \loc(\nabla^{d}b)$ is dominant, this means that it must be generically finite-to-one.
Hence, setting $c:=\nabla^{d+1}(b)$, we have that $c\in\acl(\nabla^d(b))$.
As in the proof of Lemma~\ref{extacl}, it follows that if $Y=\loc(c/\nabla^d(b))$ then $\tau_cY$ is a singleton.
Since $Y$ is defined over $\nabla^d(b)$, the prolongation space $\tau Y$ is defined over $\nabla^{d+1}(b)=c$, and hence also $\tau_cY$ is defined over $c$.
By Lemma~\ref{pronab}, $\nabla(c)\in\tau_cY$.
So $\nabla(c)\in\dcl(c)$.
By quantifier eliimnation in $\ccm$ we can write $\nabla(c)=v(c)$ for some definable meromorphic map $v$.
Then, setting $X:=\loc(c)$,  we have that $v:X\to TX$ is a section to the tangent space of $X$.
That is, $(X,v)$ is a meromorphic vector field and $q=\tp(c)$ is its generic type.
Finally, observe that $b$ and $c=\nabla^{d+1}(b)$ are interdefinable over the empty set, so that $p$ and $q$ are interdefinable types.
\end{proof}

The upshot is that the finite-dimensional fragment of $\dccm$, over the empty set, captures precisely the bimeromorphic geometry of meromorphic vector fields.

\begin{remark}
We have restricted our attention in this discussion to the empty set for brevity; we could have worked more generally over arbitrary parameters $A$.
The result would be that the finite-dimensional types over $A$ are precisely, up to interdefinability, the generic types of {\em meromorphic $D$-varieties} over $A$.
We leave it to the reader to both articulate precisely, and verify, this claim.
\end{remark}

\vfill
\pagebreak


\end{document}